\pdfoutput=1
\RequirePackage{amsmath}
\documentclass[letterpaper,11pt]{article}

\usepackage[margin=1in]{geometry}

\usepackage{marvosym}

\usepackage[colorlinks,citecolor=red,urlcolor=blue,bookmarks=false,hypertexnames=true]{hyperref}

\usepackage{amsfonts}
\usepackage{amsmath}
\usepackage{amsthm}

\usepackage{algorithm}
\usepackage[noend]{algpseudocode}

\usepackage{cite}
\usepackage{xcolor}
\usepackage{color}
\usepackage{graphicx}

\newcommand{\R}{\mathbb{R}} 

\newcommand{\cD}{{\cal D}}

\newcommand{\mA}{{\bf A}}

\newcommand{\mD}{{\bf D}}

\newcommand{\mI}{{\bf I}}

\newcommand{\mM}{{\bf M}}

\newcommand{\mP}{{\bf P}}
\newcommand{\mQ}{{\bf Q}}

\newcommand{\mS}{{\bf S}}

\newcommand{\mW}{{\bf W}}

\newcommand{\trans}{\top}

\newcommand{\eqdef}{\overset{\text{def}}{=}} 
\newcommand{\ve}[2]{\langle #1 ,  #2 \rangle} 
\newcommand{\bve}[2]{\bigg\langle #1 ,  #2 \bigg\rangle} 
\newcommand{\norm}[1]{\lVert#1\rVert}      
\newcommand{\abs}[1]{\left|#1\right|}

\newcommand{\nthreads}{q}

\newcommand{\Prob}[1]{\mathbb{P}(#1)}

\newcommand{\siggen}[3]{#1_{#2}\left(#3\right)}
\newcommand{\sigmin}[1]{\siggen{\sigma}{\min}{#1}}
\newcommand{\sigmax}[1]{\siggen{\sigma}{\max}{#1}}

\newcommand{\siggenA}[2]{#1_{#2}(\mA)}

\newcommand{\sigminsqA}{\siggenA{\sigma^2}{\min}}
\newcommand{\sigmaxsqA}{\siggenA{\sigma^2}{\max}}

\DeclareMathOperator*{\argmin}{arg\,min}

\newcommand{\thup}{^{\text{th}}} %
\newcommand{\ik}{{i_k}}

\newcommand{\Diag}[1]{\mathbf{Diag}\left( #1\right)}

\newcommand{\expSB}[2]{{ \mathbf{E}}_{#1}\left[#2\right] } 
\newcommand{\E}[1]{{\mathbb{E}}\left[#1\right] }

\newtheorem{theorem}{Theorem} 
\newtheorem{assumption}{Assumption} 
\newtheorem{definition}{Definition}
\newtheorem{lemma}{Lemma}
\newtheorem{corollary}{Corollary}

\usepackage[symbol]{footmisc}

\makeatletter
\let\cl@chapter\undefined
\makeatletter

\usepackage{cleveref}
\crefname{subsection}{subsection}{subsections}
\newcommand{\renewtheorem}[1]{%
  \expandafter\let\csname #1\endcsname\relax
  \expandafter\let\csname c@#1\endcsname\relax
  \expandafter\let\csname end#1\endcsname\relax
  \newtheorem{#1}%
}
\renewtheorem{assumption}{Assumption}

\usepackage{subcaption}
\begin{document}

\title{Randomized Kaczmarz with Averaging\footnotemark[1]}


\author{Jacob D. Moorman\footnotemark[2]~\textsuperscript{\Letter} \and
        Thomas K. Tu\footnotemark[2] \and 
        Denali Molitor\footnotemark[2] \and
        Deanna Needell\footnotemark[2]
}


\footnotetext{\noindent\footnotemark[1]~Molitor, Moorman, and Needell were partially supported by NSF CAREER DMS $\#1348721$ and NSF BIGDATA $\#1740325$. Moorman was additionally supported by NSF grant DGE-$1829071$. Tu was supported by DARPA under agreement number FA8750-18-2-0066.\\
\indent \indent \footnotemark[2]~Department of Mathematics, University of California, Los Angeles, Los Angeles, CA 90095-0001, USA\\
\indent \indent \textsuperscript{\Letter}~jdmoorman@math.ucla.edu 
           }


\maketitle

\begin{abstract}
The randomized Kaczmarz (RK) method is an iterative method for approximating the least-squares solution of large linear systems of equations. The standard RK method uses sequential updates, making parallel computation difficult. Here, we study a parallel version of RK where a weighted average of independent updates is used. We analyze the convergence of RK with averaging and demonstrate its performance empirically. We show that as the number of threads increases, the rate of convergence improves and the convergence horizon for inconsistent systems decreases.
\end{abstract}

\section{Introduction}
In computed tomography, image processing, machine learning, and many other fields, a common problem is that of finding solutions to large linear systems of equations.
Given $\mA\in\R^{m\times n}$ and $b\in\R^m$, we aim to find $x \in \R^n$ which solves the linear system of equations
\begin{equation}\label{eqn:system}
    \mA x = b.
\end{equation} 
We will generally assume the system is overdetermined, with $m \gg n$. For simplicity, we assume throughout that $\mA$ has full rank so that the solution is unique when it exists. However, this assumption can be relaxed by choosing the solution with least-norm when multiple solutions exist.

When a solution to \Cref{eqn:system} exists, we denote the solution by $x^\star$ and refer to the problem as \textit{consistent}. Otherwise, the problem is \textit{inconsistent}, and $x^\star$ instead denotes the \textit{least-squares} solution
\begin{equation*}
x^\star \eqdef \underset{x\in \R^n}{\argmin} \frac{1}{2}\norm{b - \mA x}_2^2.
\end{equation*}
The least-squares solution can be equivalently written as $x^\star = \mA^\dagger b,$ where $\mA^\dagger$ is the Moore-Penrose pseudoinverse of $\mA.$ We denote the least-squares \textit{residual} as $r^\star \eqdef b - \mA x^\star$, which is zero for consistent systems.

\subsection{Randomized Kaczmarz}\label{subsec:algo}
Randomized Kaczmarz (RK) is a popular iterative method for approximating the least-squares solution of large, overdetermined linear systems \cite{Kaczmarz1937, Strohmer2009}. At each iteration, an equation is chosen at random from the system in \Cref{eqn:system} and the current iterate is projected onto the solution space of that equation. In a relaxed variant of RK, a step is taken in the direction of this projection with the size of the step depending on a relaxation parameter.

Let $x^k$ be the $k\thup$ iterate. We use $\mA_i$ to denote the $i\thup$ row of $\mA$ and $\norm{\cdot} \eqdef \norm{\cdot}_2$. The \textit{relaxed RK} update is given by 
\begin{equation}\label{eqn:kaczmarz-update}
    x^{k+1} = x^{k} - \lambda_{k,\ik}\frac{\mA_\ik x^{k} - b_\ik}{\norm{\mA_\ik}^2} \mA_\ik^\trans,
\end{equation}
where $\ik$ is sampled from some fixed distribution $\cD$ at each iteration and $\lambda_{k,i}$ are relaxation parameters \cite{CaiRelaxedRK}. Fixing the relaxation parameters $\lambda_{k,i} = 1$ for all iterations $k$ and indices $i$ leads to the standard RK method in which one projects the current iterate $x^k$ onto the solution space of the chosen equation $\mA_{i_k} x = b_{i_k}$ at each iteration \cite{Strohmer2009}. Choosing relaxation parameters $\lambda_{k,i} \ne 1$ can be used to accelerate convergence or dampen the effect of noise in the linear system \cite{CaiRelaxedRK, hanke1990acceleration, hanke1990use}. 

For consistent systems, RK converges exponentially in expectation to the solution $x^\star$ \cite{Strohmer2009}, which when multiple solutions exist is the least-norm solution \cite{zouzias2013REK,Ma2015}. For inconsistent systems, there exists at least one equation $\mA_j x = b_j$ that is not satisfied by $x^\star$. As a result RK cannot converge for inconsistent systems, since it will occasionally project onto the solution space of such an equation. One can, however, guarantee exponential convergence in expectation to within a radius of the least-squares solution \cite{Needell2010, zouzias2013REK, Needell2012Block}. This radius is commonly referred to as the \textit{convergence horizon}.

\subsection{Randomized Kaczmarz with Averaging}

In order to take advantage of parallel computation and speed up the convergence of RK, we consider a simple extension of the RK method, where at each iteration multiple independent updates are computed in parallel and a weighted average of the updates is used.
Specifically, we write the averaged RK update
\begin{equation}\label{eqn:averaged-kaczmarz-update}
    x^{k+1} = x^{k} - \frac{1}{\nthreads}\sum_{i\in\tau_k} w_i \frac{\mA_i x^{k} - b_i}{\norm{\mA_i}^2} \mA_i^\trans,
\end{equation}
where $\tau_k$ is a random set of $\nthreads$ row indices sampled \emph{with replacement} and $w_i$ represents the weight corresponding to the $i\thup$ row.  
RK with averaging is detailed in \Cref{algo:AvgK}. If $\tau_k$ is a set of size one and the weights are chosen as $w_i = 1$ for $i=1,\ldots, m$, we recover the standard RK method. 

\begin{algorithm}
\begin{algorithmic}[1]
\State{\textbf{Input} $\mA \in \R^{m\times n}$, $b\in \R^m$, $x^0 \in \R^n$, weights $w\in \R^m$, number of maximum number of iterations $K$, distribution $\cD$, number of threads $\nthreads$}
\For{$k = 0,\dots, K-1$}
    \State $\tau_k \gets \nthreads$ indices sampled from $\cD$ 
    \State Compute $\delta \gets \frac{1}{\nthreads}\sum_{i\in\tau_k} w_i \frac{\mA_i x^{k} - b_i}{\norm{\mA_i}^2} \mA_i^\trans$ in parallel
    \State Update $x^{k+1} \gets x^{k} -\delta$
\EndFor
\State{\textbf{Output} $x^K$}
\end{algorithmic}
\caption{Randomized Kaczmarz with Averaging}\label{algo:AvgK}
\end{algorithm}

\subsection{Contributions}

We derive a general convergence result for RK with averaging, and identify the conditions required for convergence to the least-squares solution. These conditions guide the choices of weights and probabilities of row selection, up to a relaxation parameter $\alpha$.
When $\nthreads=1$ and appropriate weights and probabilities are chosen, we recover the standard convergence for RK. \cite{Strohmer2009,Needell2010,zouzias2013REK}.

For uniform weights and consistent systems, we relate RK with averaging to a more general parallel sketch-and-project method \cite{2017RichtarikStoch}. We also provide an estimate of the optimal choice for the relaxation parameter $\alpha$, and compare to the estimated optimal relaxation parameter for the sketch-and-project method \cite{2017RichtarikStoch}. Through experiments, we show that our estimate lies closer to the observed result.

\subsection{Organization}

In \Cref{sec:conv}, we analyze the convergence of RK with averaging, and state our general convergence result in \Cref{subsec:gen_thm}. In \Cref{sec:unif_weights}, we discuss the special case where the weights are chosen to be uniform, and in \Cref{sec:consistent_systems}, we discuss the special case where the system is consistent. 
In \Cref{sec:alpha_convergence}, we derive an estimate of the optimal relaxation parameter for consistent systems. In \Cref{sec:experiments}, we experimentally explore the effects of the number of threads $\nthreads$, the relaxation parameter $\alpha$, the weights $w_i$, and the distribution $\cD$ on the convergence properties of RK with averaging. 

\subsection{Related Work}

The Kaczmarz algorithm was originally proposed by Kaczmarz in 1937 \cite{Kaczmarz1937}, though it was later independently developed by researchers in computed tomography as the Algebraic Reconstruction Technique \cite{gordon1970algebraic, byrne2007applied}. The original Kaczmarz method cycles through rows in a fixed order; however, this is known to perform poorly for certain orders of the rows \cite{hamaker1978angles}. Other Kaczmarz variants \cite{xu2002method} use deterministic methods to choose the rows, but their analysis is complicated and convergence results are somewhat unintuitive.

Some randomized control methods were proposed \cite{herman1993algebraic}, but with no explicit proofs of convergence until Strohmer and Vershynin's 2009 paper \cite{Strohmer2009}, which proved RK converges linearly in expectation, with a rate directly related to geometric properties of the matrix $A$. This proof was later extended to inconsistent systems \cite{Needell2010}, showing convergence within a convergence horizon of the least-squares solution.

RK is a well-studied method with many variants. We do not provide an exhaustive review of the related literature \cite{leventhal2010randomized, zouzias2013REK, needell2013two, chen2012almost, eldar2011acceleration}, but instead only remark on some closely related parallel extensions of RK. 

Block Kaczmarz \cite{Elfving1980Block, eggermont1981block, aharoni1989block, Needell2012Block, xu2015block} randomly selects a block of rows from $\mA$ at each iteration and computes its Moore-Penrose pseudoinverse. The pseudoinverse is then applied to the relevant portion of the current residual and added to the estimate, solving the least-squares problem only on the selected block of rows. Computing the pseudoinverse, however, is costly and difficult to parallelize.

The CARP algorithm \cite{Gordon2005:CARP} also distributes rows of $\mA$ into blocks. However, instead of taking the pseudoinverse, the Kaczmarz method is then applied to the rows contained within each block. Multiple blocks are computed in parallel, and a component-averaging operator combines the approximations from each block. While CARP is shown to converge for consistent systems and to converge cyclically for inconsistent systems, no exponential convergence rate is given. 

AsyRK \cite{Wright:AsyncPRK} is an asynchronous parallel RK method that results from applying Hogwild! \cite{NiuRecht:Hogwild} to the least-squares objective. In AsyRK, each thread chooses a row $\mA_i$ at random and updates a random coordinate within the support of that row $\mA_i$ with a weighted RK update. AsyRK is shown to have exponential convergence, given conditions on the step size. Their analysis requires that $\mA$ is sparse, while we do not make this restriction.

RK falls under a more general class of methods often called sketch-and-project methods \cite{Gower2015}. For a linear system $\mA x = b,$ sketch-and-project methods iteratively project the current iterate onto the solution space of a sketched subsystem $\mS^\trans \mA x -\mS^\trans b.$ In particular, RK is a sketch-and-project method with $\mS^\trans = \mI_i$, where $\mI_i$ is the $i\thup$ row of the identity matrix. Other popular iterative methods such as coordinate descent can also be framed as sketch-and-project methods. In \cite{2017RichtarikStoch}, the authors discuss a more general version of \Cref{algo:AvgK} for sketch-and-project methods with averaging. Their analysis and discussion, however, focus on consistent systems and require uniform weights. We instead restrict our analysis to RK, but allow inconsistent systems and general weights $w_i$.

RK can also be interpreted as a subcase of stochastic gradient descent (SGD) \cite{robbins1951stochastic} applied to the loss function \cite{ NeedellSrebroWard2015}
\begin{equation*}
    F(x) = \sum_{i=1}^n f_i(x) = \sum_{i=1}^n \frac{1}{2}(\mA_i x - b_i)^2.
\end{equation*}
In this context, RK with averaging can be seen as mini-batch SGD \cite{bottou-98x,needell2016batched} with importance sampling, with the update
\begin{equation*}
    x^{k+1} = x^k - \frac{1}{\nthreads}\sum_{i\in\tau_k} \frac{w_i}{L_i} \nabla f_i(x),
\end{equation*}
where $L_i=\norm{\mA_i}^2$ is the Lipschitz constant of $\nabla f_i(x)=(\mA_i x - b_i)\mA_i^\trans$.

\section{Convergence of RK with Averaging}\label{sec:conv}
For inconsistent systems, RK satisfies the error bound
\begin{equation}\label{eqn:NoisyBound}
    \begin{split}
        &\E{\norm{e^{k+1}}^2} \le \left(1 - \frac{\sigminsqA}{\norm{\mA}_F^2}\right) \E{\norm{e^k }^2} + \frac{\norm{r^\star}^2}{\norm{\mA}_F^2}, 
    \end{split}
\end{equation}
where $e^k\eqdef x^k - x^\star$ is the error of the $k\thup$ iterate, $\sigmin{\mA}$ is the smallest nonzero singular value of $\mA$, $\norm{\mA}_F^2 = \sum_{i,j} \mA_{ij}^2$ and $r^\star$ is the least-squares residual \cite{Needell2010, zouzias2013REK}. Iterating this error bound yields 
\begin{equation*}
    \begin{split}
        &\E{\norm{e^k}^2} \le \left(1 - \frac{\sigminsqA}{\norm{\mA}_F^2}\right)^k \norm{e^0 }^2 + \frac{\norm{r^\star}^2}{\sigminsqA}. 
    \end{split}
\end{equation*}
For consistent systems the least-squares residual is $r^\star=0$ and this bound guarantees exponential convergence in expectation at a rate $1 - \frac{\sigma^2_{\min} (\mA)}{\norm{\mA}_F^2}$ \cite{Strohmer2009}. For inconsistent systems, this bound only guarantees exponential convergence in expectation to within a convergence horizon $\norm{r^\star}^2/\sigminsqA$.

 We derive a convergence result for \Cref{algo:AvgK} which is similar to \Cref{eqn:NoisyBound} and leads to a better convergence rate and a smaller convergence horizon for inconsistent systems when using uniform weights. To analyze the convergence, we begin by finding the update to the error at each iteration. Subtracting the exact solution $x^\star$ from both sides of the update rule in \Cref{eqn:kaczmarz-update} and using the fact that $\mA_i e^k - r_i^\star = \mA_i x^k - b_i$, we arrive at the error update  
\begin{equation}\label{eqn:error-update}
    e^{k+1} = e^{k} - \frac{1}{\nthreads}\sum_{i\in\tau_k} w_i \frac{\mA_i e^{k} - r^\star_i}{\norm{\mA_i}^2} \mA_i^\trans .
\end{equation}
To simplify notation, we define the following matrices.
\begin{definition}\label{defn:Mk}
Define the weighted sampling matrix
\begin{equation*}
\mM_k \eqdef \frac{1}{\nthreads}\sum_{i\in\tau_k} w_i \frac{\mI_i^\trans\mI_i}{\norm{\mA_i}^2},
\end{equation*}
where $\tau_k$ is a set of indices sampled independently from $\cD$ with replacement and $\mI$ is the identity matrix.
\end{definition}

Using \Cref{defn:Mk}, the error update from \Cref{eqn:error-update} can be rewritten as 
\begin{equation}\label{eqn:error-update-matrices}
    e^{k+1} 
    = (\mI - \mA^\trans \mM_k \mA)e^{k} + \mA^\trans \mM_k r^\star.
\end{equation}
\begin{definition}\label{dfn:PWD-mats}
Let $\Diag{ d_1, d_2, \ldots, d_m}$ denote the diagonal matrix with $d_1, d_2, \ldots d_m$ on the diagonal. 
Define the normalization matrix 
\begin{equation*}
    \mD \eqdef \Diag{\norm{\mA_1}, \norm{\mA_2}, \ldots, \norm{\mA_m}}
\end{equation*}
so that the matrix $\mD^{-1}\mA$ has rows with unit norm,
the probability matrix 
\begin{equation*}
    \mP \eqdef \Diag{p_1, p_2,\ldots, p_m},
\end{equation*}
where $p_j = \Prob{i=j}$ with $i\sim \cD,$ and
the weight matrix
\begin{equation*}
    \mW \eqdef \Diag{w_1, w_2, \ldots, w_m}.
\end{equation*}
\end{definition}

The convergence analysis additionally relies on the expectations given in
\Cref{lem:EMk}, whose proof can be found in \Cref{sec:exp_lemma_proof}.
\begin{lemma}\label{lem:EMk} 
     Let $\mM_k, \mP,\mW$, and $\mD$ be defined as in \Cref{defn:Mk,dfn:PWD-mats}. Then 
\begin{equation*}
\E{\mM_k} = \mP\mW\mD^{-2}
\end{equation*}
and 
\begin{equation*}
\begin{split}
    &\E{\mM_k^\trans \mA \mA^\trans \mM_k} 
    =  \frac{1}{\nthreads}\mP \mW^2 \mD^{-2} + \left(1 - \frac{1}{\nthreads}\right) \mP \mW\mD^{-2} \mA\mA^\trans \mP \mW \mD^{-2}.
\end{split}
\end{equation*}
\end{lemma} 

\subsection{Coupling of Weights and Probabilities}\label{subsec:coupling_weights_probs}
Note that the weighted sampling matrix $\mM_k$ is a sample average, with the number of samples being the number of threads $\nthreads$. Thus, as the number of threads $\nthreads$ goes to infinity, we have
\begin{equation*}
    \mM_k \overset{\nthreads \to \infty}{\longrightarrow} \expSB{i \sim \cD}{w_i \frac{\mI_i^\trans\mI_i}{\norm{\mA_i}^2} }=\mP\mW\mD^{-2}.
\end{equation*}
Therefore, as we take more and more threads, the averaged RK update of \Cref{eqn:averaged-kaczmarz-update} approaches the deterministic update 
\begin{equation*}
x^{k+1} = (\mI - \mA^\trans \mP\mW\mD^{-2} \mA)x^{k} + \mA^\trans \mP\mW\mD^{-2} b.
\end{equation*}
and likewise the corresponding error update in \Cref{eqn:error-update-matrices} approaches the deterministic update 
\begin{equation*}
e^{k+1} = (\mI - \mA^\trans \mP\mW\mD^{-2} \mA)e^{k} + \mA^\trans \mP\mW\mD^{-2} r^\star.
\end{equation*}
Since we want the error of the limiting averaged RK method to converge to zero, we should require that this limiting error update have the zero vector as a fixed point. Thus, we ask that
\begin{equation*}
0 = \mA^\trans \mP\mW\mD^{-2} r^\star
\end{equation*}
for any least-squares residual $r^\star$. This is guaranteed if \Cref{ass:propI} holds.

\begin{assumption}\label{ass:propI} The probability matrix $\mP$ and weight matrix $\mW$ are chosen to satisfy
\begin{equation*}
    \mP\mW\mD^{-2} = \alpha \mI.
\end{equation*}
for some scalar relaxation parameter $\alpha>0$.
\end{assumption}

\subsection{General Result}\label{subsec:gen_thm}
We now state a general convergence result for RK with averaging in \Cref{thm:expected_error}. The proof is given in \Cref{sec:expected_error_proof}. \Cref{thm:expected_error} in its general form is difficult to interpret, so we defer a detailed analysis to \Cref{sec:unif_weights} in which the assumption of uniform weights simplifies the bound significantly.
\begin{theorem}\label{thm:expected_error}
Suppose $\mP$ and $\mW$ of \Cref{dfn:PWD-mats} are chosen such that $\mP \mW \mD^{-2}= \frac{\alpha}{\norm{\mA}_F^2} \mI$ for relaxation parameter $\alpha > 0$. Then the error at each iteration of \Cref{algo:AvgK} satisfies
\begin{align*}
    &\E{\norm{e^{k+1}}^2} \le \sigma_{\max}\left(\left(\mI - \alpha \frac{\mA^\trans\mA}{\norm{\mA}_F^2}\right)^2 \right. - \left. \frac{\alpha^2}{\nthreads} \left(\frac{\mA^\trans\mA}{\norm{\mA}_F^2}\right)^2\right)\norm{e^{k}}^2 + \frac{\alpha}{\nthreads}\frac{ \norm{r^k}_\mW^2}{\norm{\mA}_F^2},
\end{align*}
where  $r^k \eqdef b - \mA x^k$ is the residual of the $k\thup$ iterate, $\norm{\cdot}_\mW^2 = \ve{\cdot}{ \mW \cdot}$ and $\norm{\mA}_F^2 = \sum_{i,j} \mA_{ij}^2$.
\end{theorem}
\noindent
Here, and for the remainder of the paper, we take the expectation $\E{\norm{e^{k+1}}^2}$ conditioned on $e_k$.

As we shall see in \Cref{sec:unif_weights}, the relaxation parameter $\alpha$ and number of threads $\nthreads$ are closely tied to both the convergence horizon and convergence rate. The convergence horizon is proportional to $\frac{\alpha^2}{\nthreads}$, so smaller $\alpha$ and larger $\nthreads$ lead to a smaller convergence horizon. Increasing the value of $\alpha$ improves the convergence rate of the algorithm up to a critical point beyond which further increasing $\alpha$ leads to slower convergence rates. Increasing the number of threads $\nthreads$ improves the convergence rate, asymptotically approaching an optimal rate as $\nthreads \to \infty$. 

\section{Uniform Weights}\label{sec:unif_weights}
We can simplify the analysis significantly if we assume that $\mW = \alpha \mI,$ or equivalently that the weights are uniform. In this case, the update for each iteration becomes 
\begin{equation*}
    x^{k+1} = x^{k} - \frac{\alpha}{\nthreads} \sum_{i\in\tau_k}  \frac{\mA_i x^{k} - b_i}{\norm{\mA_i}^2} \mA_i^\trans,
\end{equation*}
where $i\in \tau_k$ are independent samples from $\cD$ with $p_i = \frac{\norm{\mA_i}^2}{\norm{\mA}_F^2}$. 
Under these conditions, the expected error bound of \Cref{thm:expected_error} can be simplified to remove the dependence on $r^k$. This simplification leads to the more interpretable error bound given in \Cref{cor:unif_weight_err}. In particular, increasing $\nthreads$ leads to both a faster convergence rate and smaller convergence horizon. 
If the relaxation parameter $\alpha$ is chosen to be one and a single row is selected at each iteration, we arrive at the RK method \cite{Strohmer2009}. Using a relaxation parameter $\alpha$ other than one results in the relaxed RK method \cite{hanke1990use,hanke1990acceleration}.
\begin{corollary}\label{cor:unif_weight_err}
Suppose $p_i = \frac{\norm{\mA_i}^2}{\norm{\mA}_F^2}$ and $\mW = \alpha\mI$.
Then the expected error at each iteration of \Cref{algo:AvgK} satisfies
\begin{align*}
    &\E{\norm{e^{k+1}}^2} 
    \le \sigma_{\max}\left(\left(\mI - \alpha \frac{\mA^\trans\mA}{\norm{\mA}_F^2}\right)^2 
    + \frac{\alpha^2 }{\nthreads}\left( \mI -  \frac{\mA^\trans\mA}{\norm{\mA}_F^2}\right)\frac{\mA^\trans\mA}{\norm{\mA}_F^2}\right)\norm{e^{k}}^2  + \frac{\alpha^2\norm{r^\star}^2}{\nthreads\norm{\mA}_F^2} .
\end{align*}
\end{corollary}
The proof of \Cref{cor:unif_weight_err} follows immediately from \Cref{thm:expected_error} and can be found in \Cref{subsec:unif_weight_cor_proof}.

\subsubsection{Randomized Kaczmarz}\label{subsubsec:REK}
If a single row is chosen at each iteration, with $\mW = \mI$ and $p_i = \frac{\norm{\mA_i}^2}{\norm{\mA}_F^2},$ then \Cref{algo:AvgK} becomes the version of RK stated in \cite{Strohmer2009}. In this case, 
\begin{equation}
    \norm{r^k}_\mW^2 = \norm{\mA e^k}^2 + \norm{r^\star}^2.
\end{equation} Applying \Cref{thm:expected_error} leads to the following corollary, which recovers the error bound in \Cref{eqn:NoisyBound}.

\begin{corollary}\label{cor:RK_conv}
Suppose $\nthreads = 1$, $\mW=\mI$ and $p_i = \frac{\norm{\mA_i}^2}{\norm{\mA}_F^2}$. 
Then the expected error at each iteration of \Cref{algo:AvgK} satisfies
\begin{align*}
    \E{\norm{e^{k+1}}^2} 
    &\le \sigmax{\mI -  \frac{\mA^\trans\mA}{\norm{\mA}_F^2}}\norm{e^{k}}^2  + \frac{ \norm{r^\star}^2}{\norm{\mA}_F^2}
     \\ &= \left(1 -  \frac{\sigminsqA}{\norm{\mA}_F^2}\right)\norm{e^{k}}^2  + \frac{ \norm{r^\star}^2}{\norm{\mA}_F^2}.
\end{align*}
\end{corollary}
A proof of \Cref{cor:RK_conv} is included in \Cref{subsec:RK_conv_proof}.

\section{Consistent Systems}\label{sec:consistent_systems}
For consistent systems, \Cref{algo:AvgK} converges to the solution $x^\star$ exponentially in expectation with the following guaranteed convergence rate. 
\begin{corollary}\label{cor:expected_error_consist}
Suppose $\mP$ and $\mW$ of \Cref{dfn:PWD-mats} are chosen such that $\mP \mW \mD^{-2}= \frac{\alpha}{\norm{\mA}_F^2} \mI$ for some constant $\alpha > 0$. Then the error at each iteration of \Cref{algo:AvgK} satisfies
\begin{align*}
    &\E{\norm{e^{k+1}}^2} 
    \le \sigma_{\max}\left(\left(\mI - \alpha \frac{\mA^\trans\mA}{\norm{\mA}_F^2}\right)^2
    +\frac{\mA^\trans}{\norm{\mA}_F} \left(\frac{\alpha}{\nthreads} \mW - \frac{\alpha^2}{\nthreads} \frac{\mA\mA^\trans}{\norm{\mA}_F^2}\right) \frac{\mA}{\norm{\mA}_F}\right)\norm{e^{k}}^2.
\end{align*}
\end{corollary}
\Cref{cor:expected_error_consist} can be derived from the proof of \Cref{thm:expected_error} with $r^\star =0$.

\section{Suggested Relaxation Parameter $\alpha$ for Consistent Systems With Uniform Weights} \label{sec:alpha_convergence}

For consistent systems and using uniform weights, \Cref{algo:AvgK} becomes a subcase of the parallel sketch-and-project method described by {Richt{\'a}rik} and {Tak{\'a}{\v{c}}} \cite{2017RichtarikStoch}.
They suggest a choice for the relaxation parameter
\begin{equation}\label{eqn:alpha-opt-rich-takac}
\alpha = \frac{\nthreads}{1 + \left(\nthreads-1\right)\frac{\sigmaxsqA}{\norm{\mA}_F^2}}
\end{equation}
chosen to optimize their convergence rate guarantee.

Analogously, for uniform weights, we can calculate the value of $\alpha$ to minimize the bound given in \Cref{cor:unif_weight_err}. 

\begin{theorem}\label{thm:unif_opt_alpha}
Suppose $p_i = \frac{\norm{\mA_i}^2}{\norm{\mA}_F^2}$ and $\mW = \alpha\mI$. Then, the relaxation parameter $\alpha$ which yields the fastest convergence rate guarantee in \Cref{cor:unif_weight_err} is
\begin{equation*}\label{eqn:alpha-opt}
    \alpha^\star = 
     \begin{cases}
       \frac{\nthreads}{1+(\nthreads-1)s_{\min}},    &1-(\nthreads-1)(s_{\max}-s_{\min})\geq 0,\\
       \frac{2 \nthreads}{1 + \left(\nthreads-1\right) \left(s_{\min} + s_{\max}\right)}, & 1-(\nthreads-1)(s_{\max}-s_{\min})<0
     \end{cases}
\end{equation*}
where $s_{\min}=\frac{\sigminsqA}{\norm{\mA}_F^2}$ and $s_{\max}=\frac{\sigmaxsqA}{\norm{\mA}_F^2}$.
\end{theorem}
The proof of this result can be found in \Cref{sec:unif_opt_alpha_proof}.

When $\nthreads=1$, the second condition cannot hold, and so only the first formula is used. Plugging in $\nthreads=1$, the term that depends on $s_min$ vanishes and we get that $\alpha^\star=1$. When $\nthreads>1$, we can divide by $\nthreads-1$ and express the condition in terms of the spectral gap as $s_{\max}-s_{\min} \leq \frac{1}{\nthreads-1}$. For matrices where the spectral gap is positive, we can also view this as a condition on the number of threads, $q \leq 1+\frac{1}{s_{\max}-s_{\min}}$. We see that for low numbers of threads, the first form is used, while for high numbers of threads, the second is used.

Note that this differs from the relaxation parameter $\alpha$ suggested by {Richt{\'a}rik} and {Tak{\'a}{\v{c}}} \cite{2017RichtarikStoch}, given in \Cref{eqn:alpha-opt-rich-takac}.
This is due to the fact that our convergence rate guarantee is tighter, and thus we expect that our suggested relaxation parameter $\alpha$ should be closer to the truly optimal value. 
We compare these two choices of the relaxation parameter $\alpha$ experimentally in \Cref{subsec:effect_of_alpha} and show that our suggested relaxation parameter $\alpha^\star$ is indeed closer to the true optimal value, especially for large numbers of threads $\nthreads$.

\section{Experiments}\label{sec:experiments}
We present several experiments to demonstrate the convergence of \Cref{algo:AvgK} under various conditions. In particular, we study the effects of the number of threads $\nthreads$, the relaxation parameter $\alpha$, the weight matrix $\mW$, and the probability matrix $\mP$.

\subsection{Procedure}
For each experiment, we run $100$ independent trials each starting with the initial iterate $x^0=0$ and average the squared error norms $\norm{e^k}^2$ across the trials. We sample $\mA$ from $100 \times 10$ standard Gaussian matrices and least-squares solution $x^\star$ from $10$-dimensional standard Gaussian vectors, normalized so that $\norm{x^\star}=1$. To form inconsistent systems, we generate the least-squares residual $r^\star$ as a Gaussian vector orthogonal to the range of $\mA$, also normalized so that $\norm{r^\star}=1$. Finally, $b$ is computed as $r^\star + \mA x^\star$. 

\subsection{The Effect of the Number of Threads}
\begin{figure}
  \begin{subfigure}[t]{0.45\textwidth}
    \includegraphics[width=\textwidth]{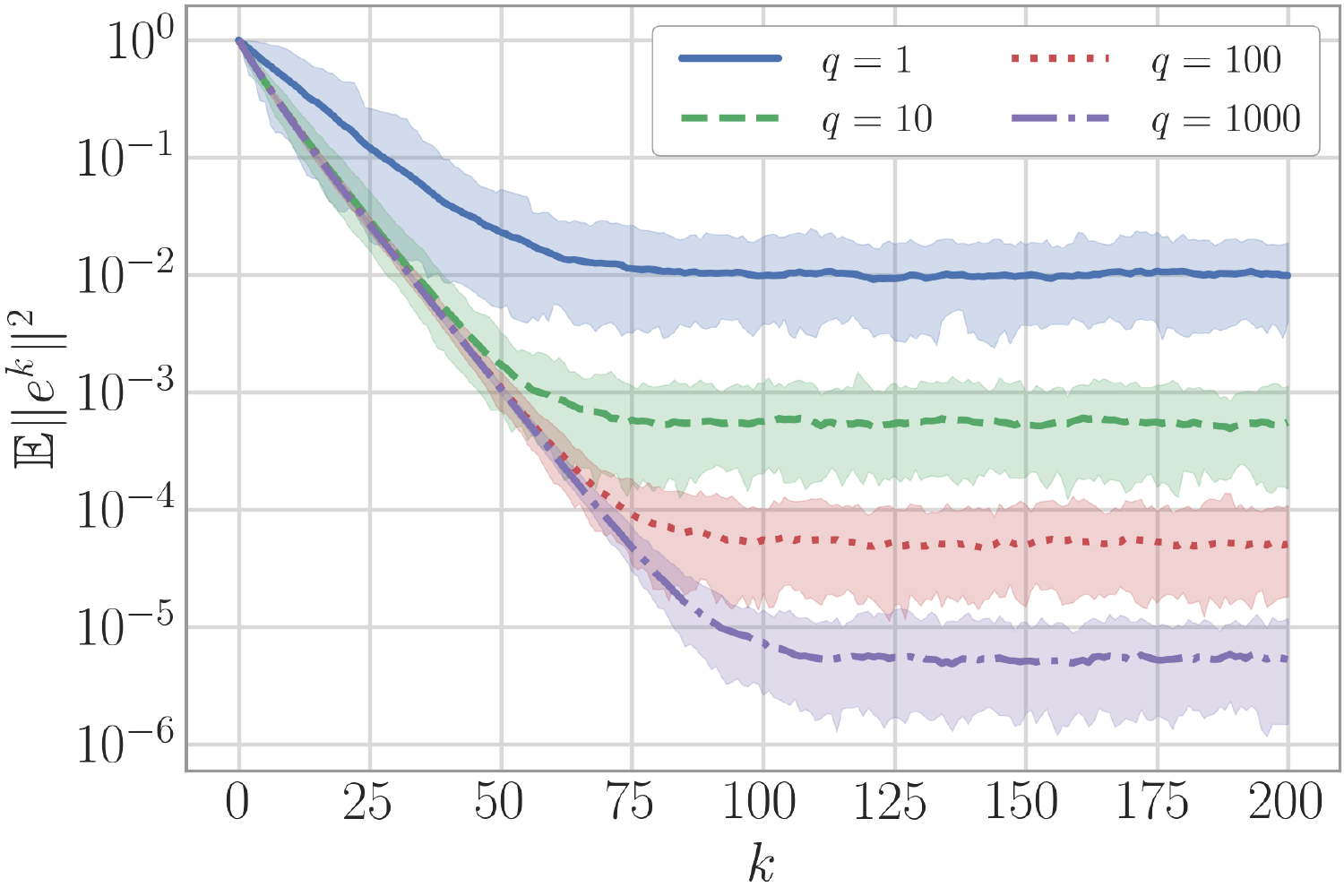}
    \vspace{-0.2 in}
    \caption{Uniform weights $w_i = 1$ and probabilities proportional to squared row norms $p_i = \frac{\norm{\mA_i}^2}{\norm{\mA}_F^2}$.}
    \label{fig:unif-weights-rownorm-probs-vary-tau}
  \end{subfigure}
  \hspace{0.05\textwidth}
  \begin{subfigure}[t]{0.45\textwidth}
    \includegraphics[width=\textwidth]{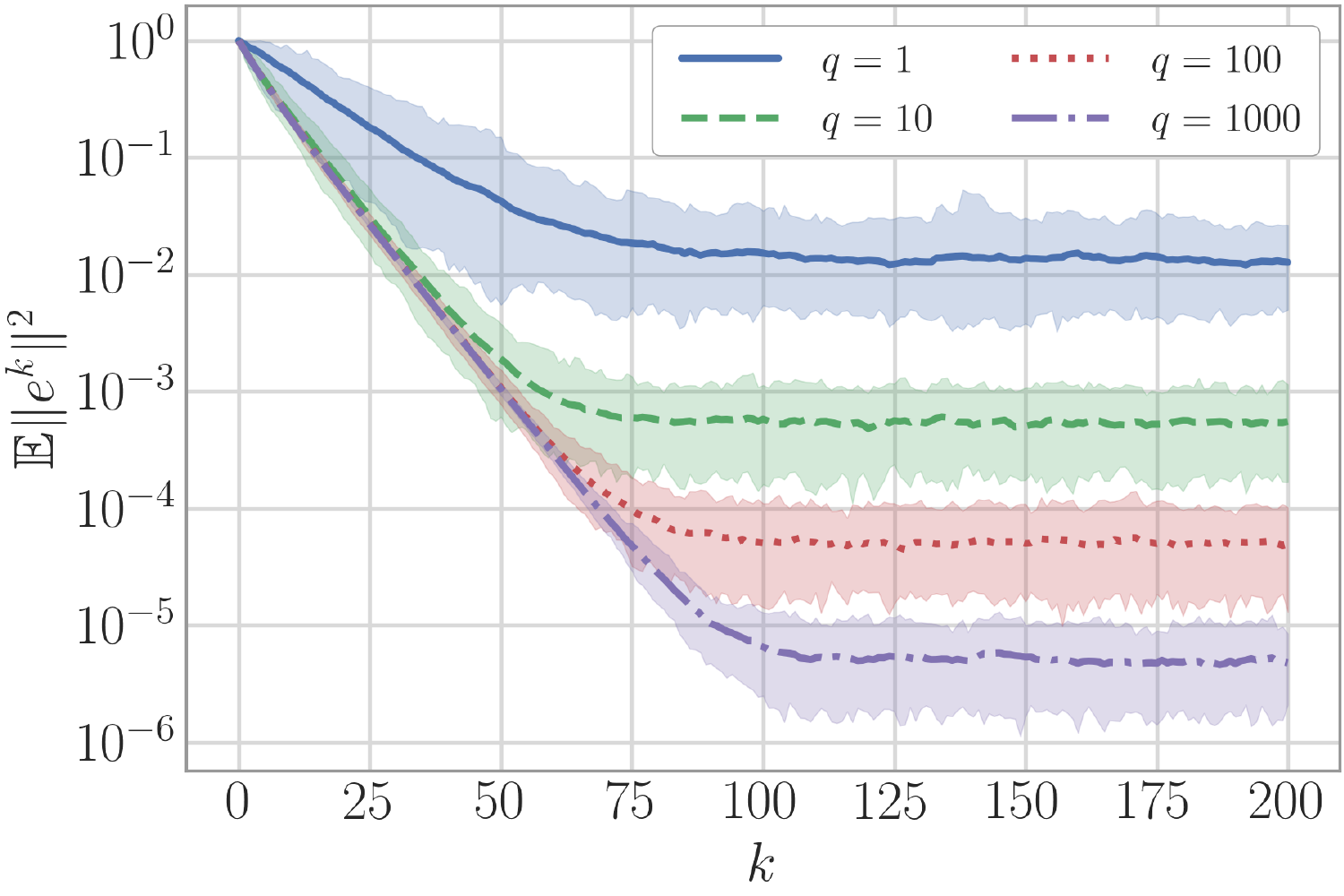}
    \vspace{-0.2 in}
    \caption{Weights proportional to squared row norms $w_i=m\frac{\norm{\mA_i}^2}{\norm{\mA}_F^2}$  and uniform probabilities $p_i = \frac{1}{m}$.}
    \label{fig:rownorm-weights-unif-probs-vary-tau}
  \end{subfigure}
  \begin{subfigure}[t]{0.45\textwidth}
    \includegraphics[width=\textwidth]{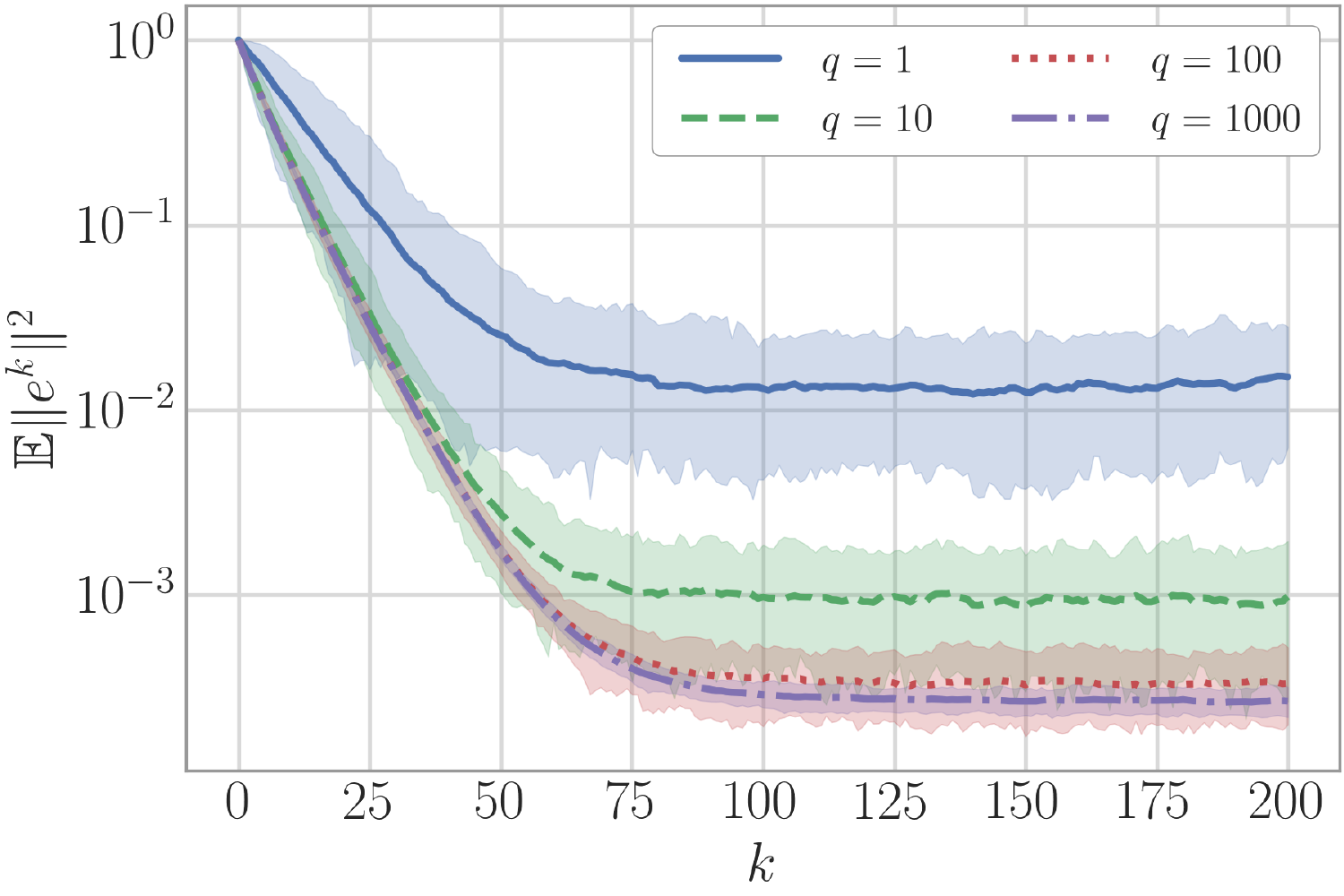}
    \vspace{-0.2 in}
    \caption{Uniform weights $w_i=1$ and uniform probabilities $p_i = \frac{1}{m}$.}
    \label{fig:unif-weights-unif-probs-vary-tau}
  \end{subfigure}
  \centering
  \caption{The effect of the number of threads on the average squared error norm vs iteration for \Cref{algo:AvgK} applied to inconsistent systems. The weights $w_i$ and probabilities $p_i$ in \subref{fig:unif-weights-rownorm-probs-vary-tau} and \subref{fig:rownorm-weights-unif-probs-vary-tau} satisfy \Cref{ass:propI},
  while in \subref{fig:unif-weights-unif-probs-vary-tau} they do not. Shaded regions are $5\thup$ and $95\thup$ percentiles, measured over 100 trials.}
  \centering
  \label{fig:vary-tau}
\end{figure}

In \Cref{fig:vary-tau}, we see the effects of the number of threads $\nthreads$ on the approximation error of \Cref{algo:AvgK} for different choices of the weight matrices $\mW$ and probability matrices $\mP$. In
\Cref{fig:unif-weights-rownorm-probs-vary-tau,fig:rownorm-weights-unif-probs-vary-tau}, $\mW$ and $\mP$ satisfy \Cref{ass:propI}, while in \Cref{fig:unif-weights-unif-probs-vary-tau} they do not.

In \Cref{fig:unif-weights-rownorm-probs-vary-tau,fig:rownorm-weights-unif-probs-vary-tau}, as the number of threads $\nthreads$ increases by a factor of ten, we see a corresponding decrease in the magnitude of the convergence horizon by approximately the same factor. This result corroborates what we expect based on \Cref{thm:expected_error} and \Cref{cor:unif_weight_err}. For \Cref{fig:unif-weights-unif-probs-vary-tau}, we do not see the same consistent decrease in the magnitude of the convergence horizon. As $\nthreads$ increases, for weight matrices $\mW$ and probability matrices $\mP$ that do not satisfy \Cref{ass:propI}, the iterates $x^k$ approach a weighted least-squares solution instead of the desired least-squares solution $x^\star$ (see \Cref{subsec:coupling_weights_probs}).

The rate of convergence in \Cref{fig:vary-tau} also improves as the number of threads $\nthreads$ increases. As $\nthreads$ increases, we see diminishing returns in the convergence rate. 
We expect this behavior based on the dependence on $\frac{1}{\nthreads}$ in \Cref{thm:expected_error} and \Cref{cor:unif_weight_err}.

\subsection{The Effect of the Relaxation Parameter $\alpha$} \label{subsec:effect_of_alpha}

\begin{figure}
  \begin{subfigure}[t]{0.45\textwidth}
    \includegraphics[width=\textwidth]{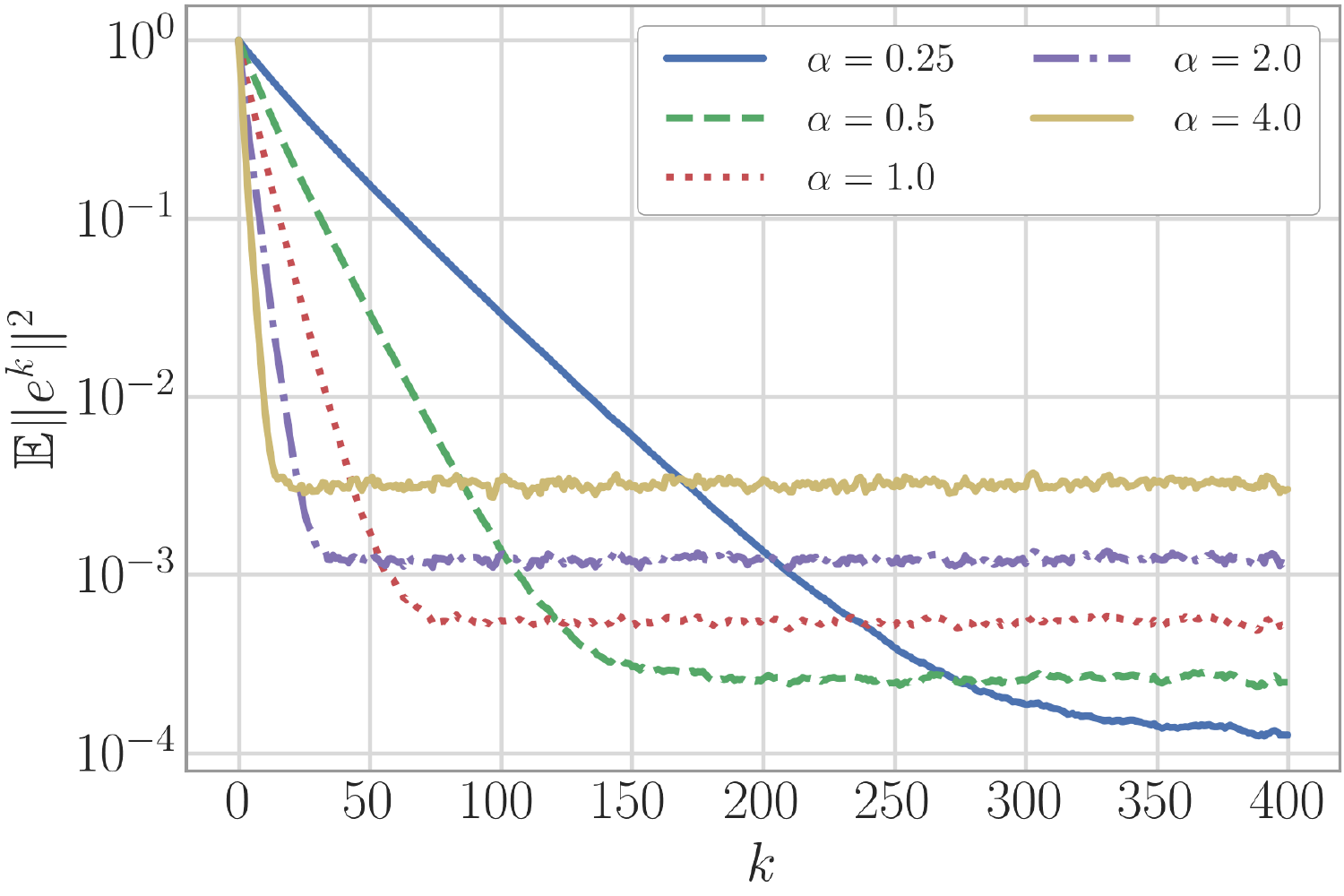}
    \vspace{-0.2 in}
    \caption{Uniform weights $w_i = \alpha$, probabilities proportional to squared row norms $p_i = \frac{\norm{\mA_i}^2}{\norm{\mA}_F^2}$, and number of threads $\nthreads=10$.}
    \label{fig:unif-weights-rownorm-probs-vary-alpha}
  \end{subfigure}
  \hspace{0.05\textwidth}
  \begin{subfigure}[t]{0.45\textwidth}
    \includegraphics[width=\textwidth]{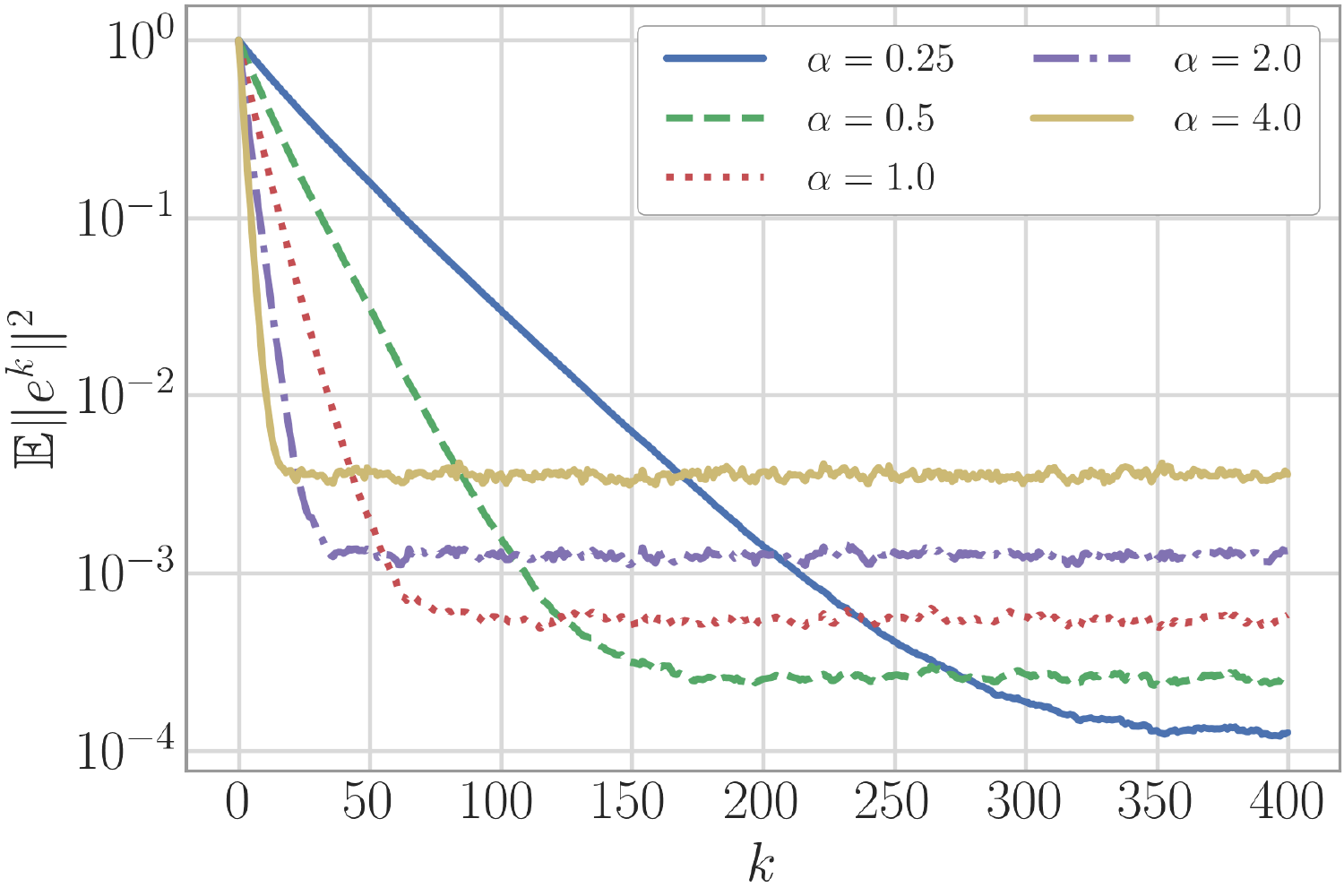}
    \vspace{-0.2 in}
    \caption{Weights proportional to squared row norms $w_i=\alpha m\frac{\norm{\mA_i}^2}{\norm{\mA}_F^2}$, uniform probabilities $p_i = \frac{1}{m}$, and number of threads $\nthreads=10$.}
    \label{fig:rownorm-weights-unif-probs-vary-alpha}
  \end{subfigure}
  \centering
  \caption{The effect of the relaxation parameter $\alpha$ on the average squared error norm vs iteration for \Cref{algo:AvgK} applied to inconsistent systems.}
  \centering
  \label{fig:coupled-vary-alpha}
\end{figure}

In \Cref{fig:coupled-vary-alpha}, we observe the effect on the convergence rate and convergence horizon as we vary the relaxation parameter $\alpha$. From \Cref{thm:expected_error}, we expect that the convergence horizon increases with $\alpha$ and indeed observe this experimentally. The squared norms of the errors behave similarly as $\alpha$ varies for both sets of weights and probabilities considered, each of which satisfy \Cref{ass:propI}. 

\begin{figure}[t]
  \begin{subfigure}[t]{0.45\textwidth}
    \includegraphics[width=\textwidth]{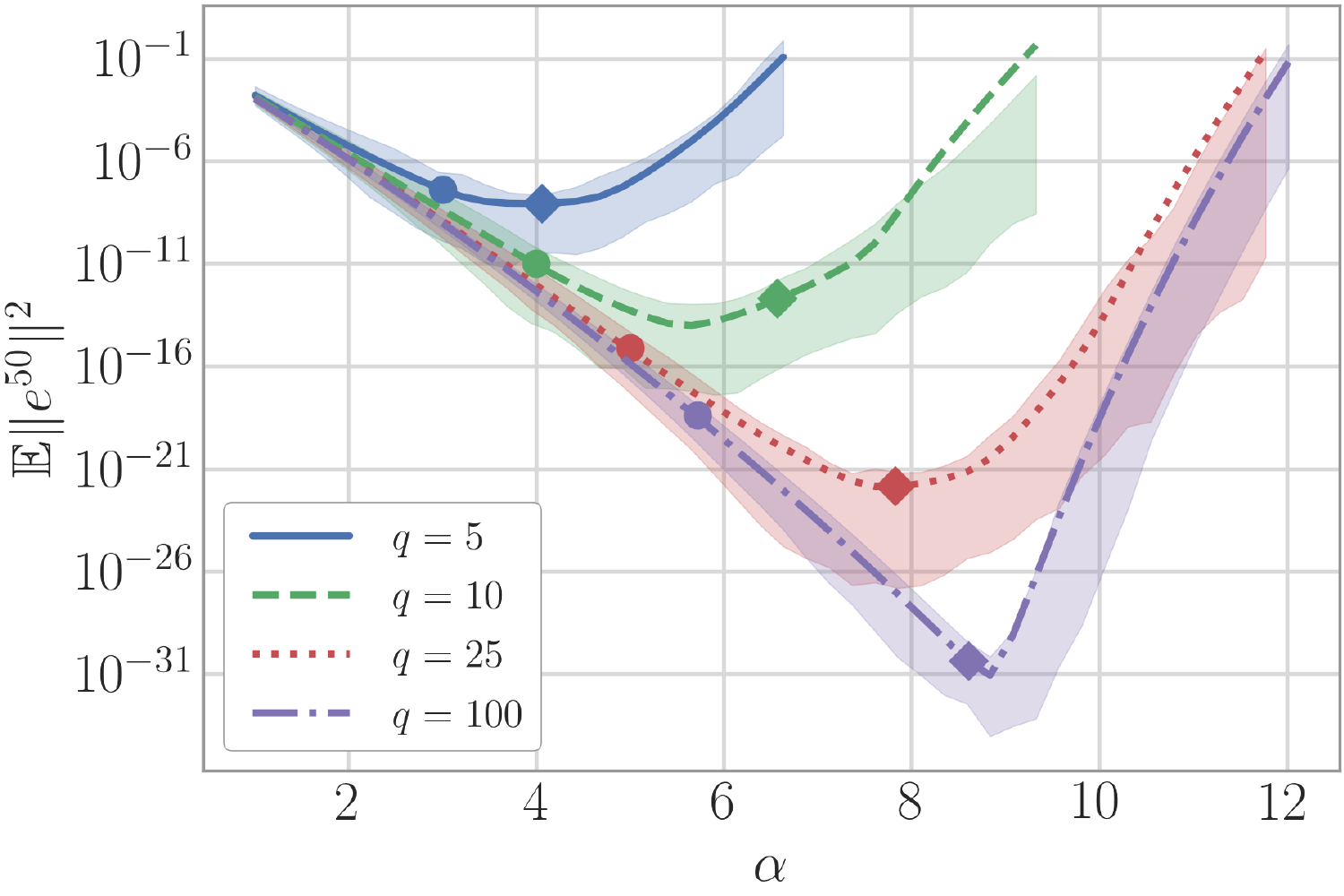}
    \vspace{-0.2 in}
    \caption{$w_i = \alpha$, $p_i = \frac{\norm{\mA_i}^2}{\norm{\mA}_F^2}$.}
    \label{fig:unif-weights-k-vs-alpha-vary-tau}
  \end{subfigure}
  \hspace{0.05\textwidth}
  \begin{subfigure}[t]{0.45\textwidth}
    \includegraphics[width=\textwidth]{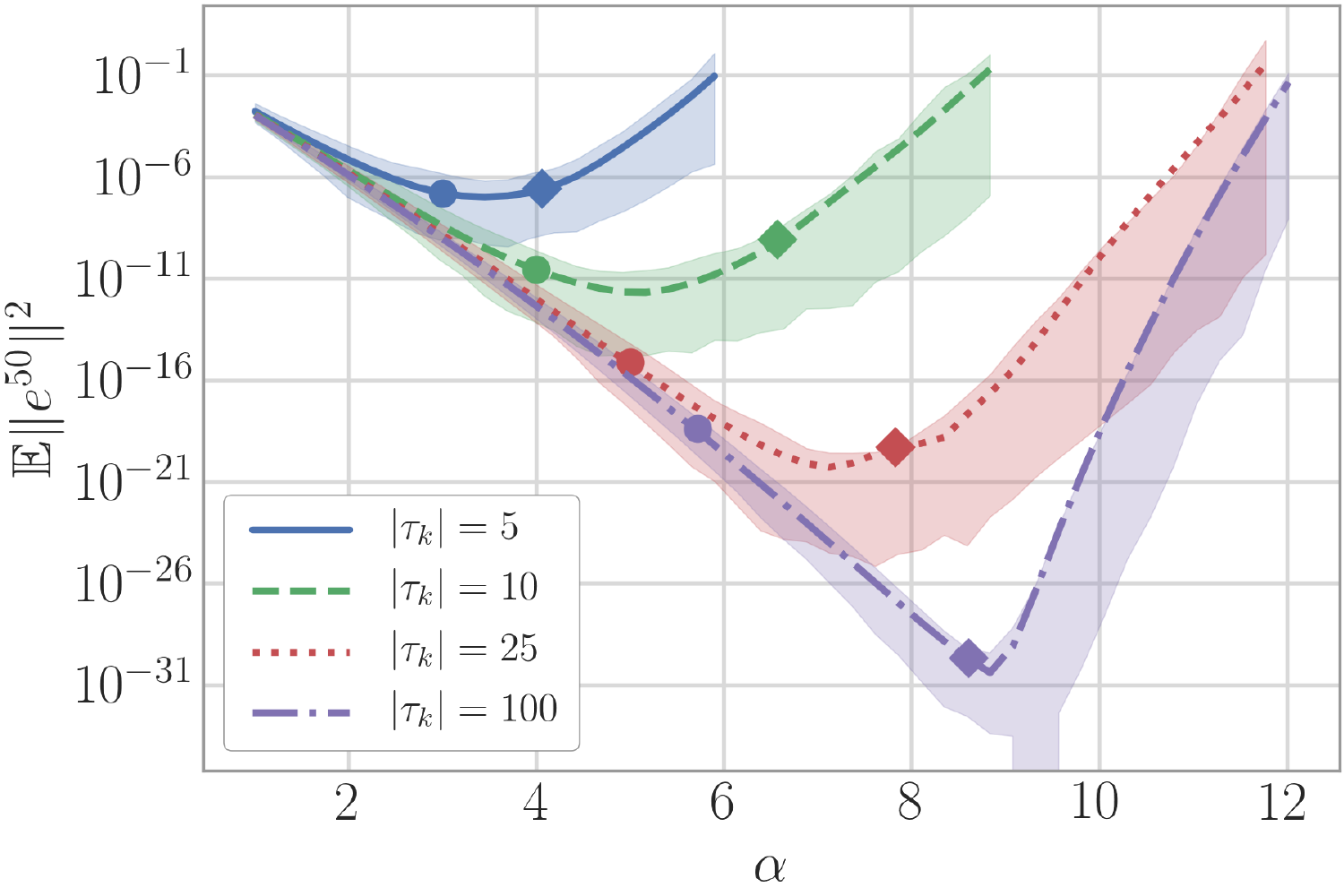}
    \vspace{-0.2 in}
    \caption{$w_i=\alpha m\frac{\norm{\mA_i}^2}{\norm{\mA}_F^2}$, $p_i = \frac{1}{m}$.}
    \label{fig:rownorm-weights-k-vs-alpha-vary-tau}
  \end{subfigure}
  \centering
  \caption{Squared error norm after 50 iterations of \Cref{algo:AvgK} on consistent systems for various choices of relaxation parameter $\alpha$. Shaded regions are the $5\thup$ and $95\thup$ percentiles, measured over 100 trials. Diamond markers are estimates of the optimal alpha using \Cref{thm:unif_opt_alpha}, and circle markers are estimates using the formula from Richt{\'a}rik and Tak{\'a}{\v{c}} \cite{2017RichtarikStoch}}
  \label{fig:errsq-vary-alpha}
\end{figure}

For larger values of the relaxation parameter $\alpha$, the convergence rate for \Cref{algo:AvgK} eventually decreases and the method can ultimately diverge. This behavior can be seen in \Cref{fig:errsq-vary-alpha}, which plots the squared error norm after 100 iterations for consistent Gaussian systems, various $\alpha$, and various numbers of threads $\nthreads$. In \Cref{fig:unif-weights-k-vs-alpha-vary-tau}, we use uniform weights $w_i=\alpha$ with probabilities proportional to the squared row norms $p_i = \frac{\norm{\mA_i}^2}{\norm{\mA}_F^2}$, and in \Cref{fig:rownorm-weights-k-vs-alpha-vary-tau}, we use weights proportional to the row norms $w_i=\alpha m\frac{\norm{\mA_i}^2}{\norm{\mA}_F^2}$ with uniform probabilities $p_i = \frac{1}{m}$.

For each value of $\nthreads$, we plot two markers on the curve to show the estimated optimal values of $\alpha$. The diamond markers are optimal values of $\alpha$ computed using \Cref{thm:unif_opt_alpha}, and the circle markers are optimal values of $\alpha$ using the formula from Richt{\'a}rik and Tak{\'a}{\v{c}} \cite{2017RichtarikStoch}. These values are also contained in \Cref{tab:alpha-opt}.
In terms of the number of iterations required, we find that the optimal value for $\alpha$ increases with $\nthreads$. 
Comparing the $\alpha$ values from \cite{2017RichtarikStoch} with the $\alpha$ that minimize the curves in \Cref{fig:errsq-vary-alpha}, we find that these values generally underestimate the optimal $\alpha$ that we observe experimentally. In comparison, the optimal $\alpha$ calculated using \Cref{thm:unif_opt_alpha} are much closer to the observed optimal values of $\alpha$, especially for high $\nthreads$

\begin{table}[h!]
    \caption{Calculated optimal $\alpha^\star$ for matrix $\mA$ used in \Cref{fig:unif-weights-k-vs-alpha-vary-tau}.}
    \label{tab:alpha-opt}
  \begin{center}
    \begin{tabular}{c|c|c|c|c}
        & $\nthreads=5$ & $\nthreads=10$ & $\nthreads=25$ & $\nthreads=100$\\
      \hline
      $\alpha$ (Eqn~\ref{eqn:alpha-opt-rich-takac}) [Richt{\'a}rik et al.] & 3.00 & 4.00 & 5.00  & 5.72\\
      \hline
      $\alpha^\star$ (our \Cref{thm:unif_opt_alpha}) & 4.06 & 6.57 & 7.83 & 8.61
    \end{tabular}
  \end{center}
\end{table}

\begin{figure}[t]
  \begin{subfigure}[t]{0.45\textwidth}
    \includegraphics[width=\textwidth]{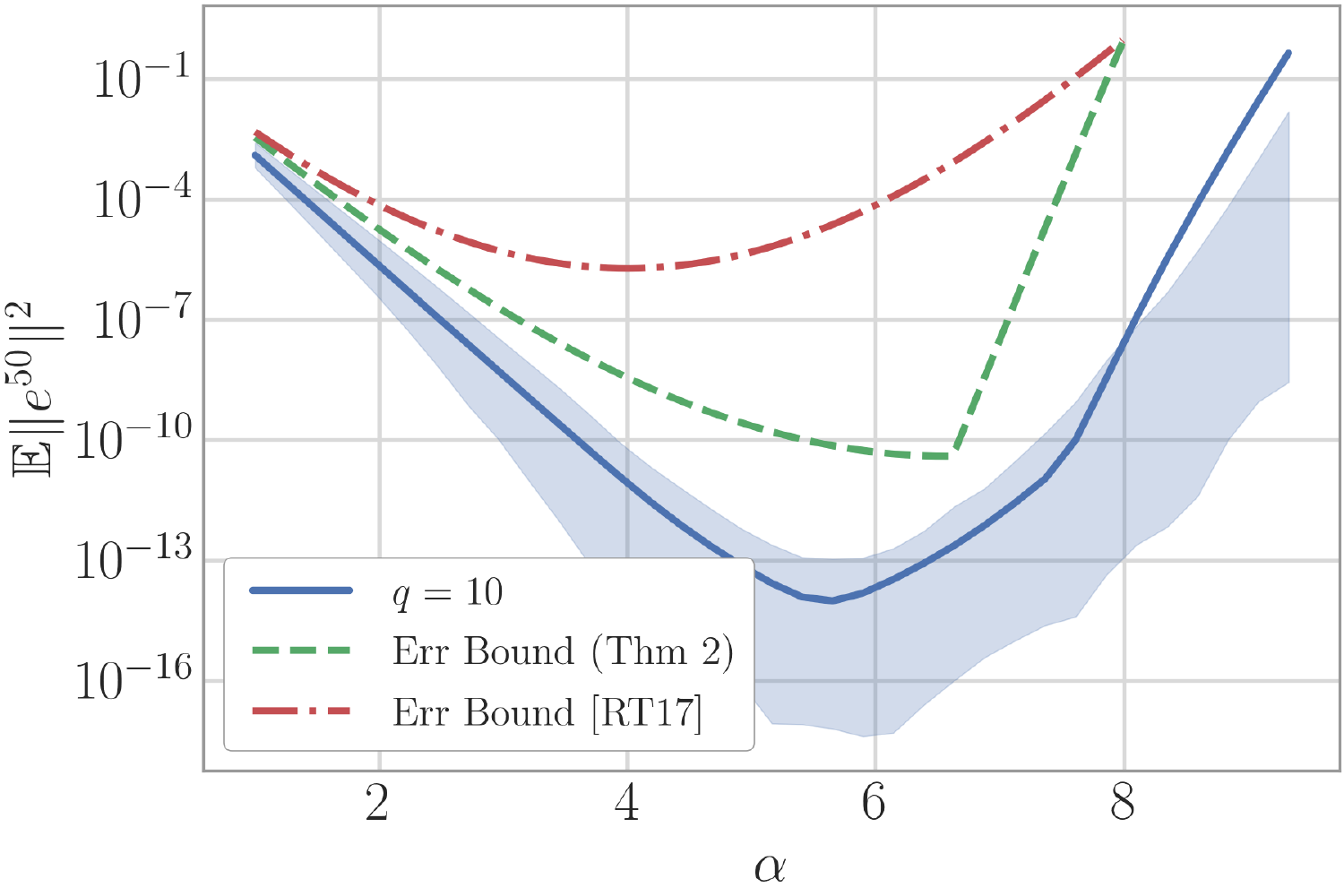}
    \vspace{-0.2 in}
    \caption{$\nthreads=10$. }
    \label{fig:err_bound_comp_tau_10}
  \end{subfigure}
  \hspace{0.05\textwidth}
  \begin{subfigure}[t]{0.45\textwidth}
    \includegraphics[width=\textwidth]{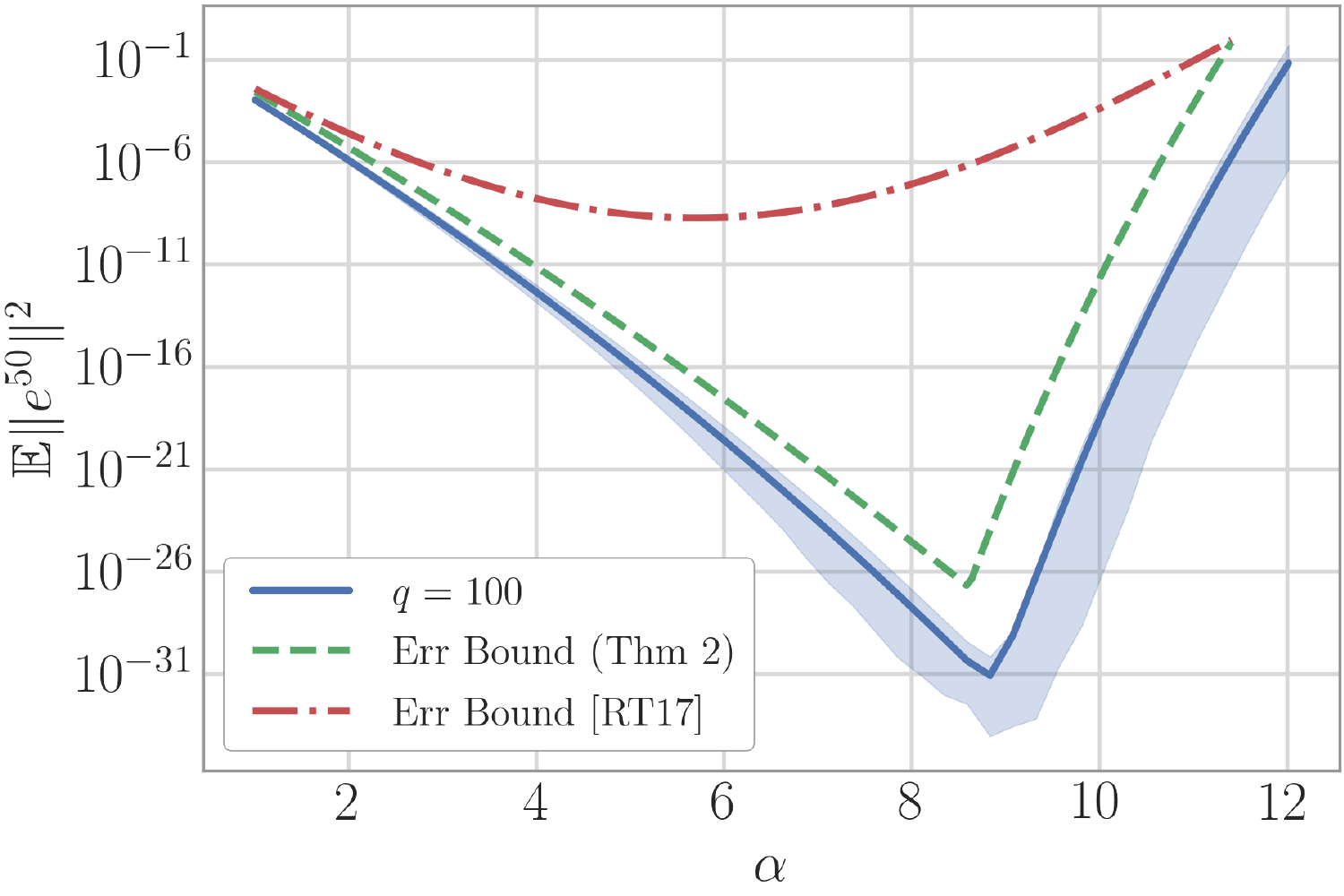}
    \vspace{-0.2 in}
    \caption{$\nthreads=100$.}
    \label{fig:err_bound_comp_tau_100}
  \end{subfigure}
  \centering
  \caption{Squared error norm after 50 iterations of \Cref{algo:AvgK} on consistent systems for various choices of relaxation parameter $\alpha$. Uniform weights $w_i = \alpha$ and probabilities proportional to squared row norms $p_i = \frac{\norm{\mA_i}^2}{\norm{\mA}_F^2}$.}
  \label{fig:err_bound_comp}
\end{figure}

We believe this is due to our bound being relatively tighter than \Cref{eqn:alpha-opt-rich-takac}. In \Cref{fig:err_bound_comp_tau_10,fig:err_bound_comp_tau_100}, we plot the error bounds produced by \Cref{eqn:alpha-opt-rich-takac} and \Cref{eqn:alpha-opt} after 50 iterations for $\nthreads=10$ and $\nthreads=100$. We observe that as the number of threads increases, our bound approaches the empirical result.

\section{Conclusion}
We prove a general error bound for RK with averaging given in \Cref{algo:AvgK} in terms of the number of threads $\nthreads$ and a relaxation parameter $\alpha$. We find a natural coupling between the probability matrix $\mP$ and the weight matrix $\mW$ that leads to a reduced convergence horizon. We demonstrate that for uniform weights, i.e.\ $\mW \propto \mI$, the rate of convergence and convergence horizon for \Cref{algo:AvgK} improve both in theory and practice as the number of threads $\nthreads$ increases. Based on the error bound, we also derive an optimal value for the relaxation parameter $\alpha$ which increases convergence speed, and compare with existing results. 

\bibliographystyle{plain}
\bibliography{aveKacz}

\appendix
\section{Proof of \Cref{lem:EMk}}\label{sec:exp_lemma_proof}
Expanding the definition of the weighted sampling matrix $\mM_k$ as a weighted average of the i.i.d. sampling matrices $\frac{\mI_i^\trans\mI_i}{\norm{\mA_i}^2}$, we see that
\begin{equation*}
\E{\mM_k} = \E{\frac{1}{\nthreads}\sum_{i\in\tau_k} w_i \frac{\mI_i^\trans\mI_i}{\norm{\mA_i}^2}} = \E{w_i \frac{\mI_i^\trans\mI_i}{\norm{\mA_i}^2}} = \sum_{i=1}^m p_i w_i \frac{\mI_i^\trans\mI_i}{\norm{\mA_i}^2} = \mP\mW\mD^{-2}.
\end{equation*}
Likewise, we can compute
\begin{align*}
    \E{\mM_k^\trans \mA \mA^\trans \mM_k} &= \E{\left(\frac{1}{\nthreads}\sum_{i\in\tau_k} w_i \frac{\mI_i^\trans\mA_i}{\norm{\mA_i}^2}\right)\left(\frac{1}{\nthreads}\sum_{j\in\tau_k} w_j \frac{\mA_j^\trans\mI_j}{\norm{\mA_j}^2}\right)} \\
    & = \frac{1}{\nthreads} \E{\left(w_i \frac{\mI_i^\trans\mA_i}{\norm{\mA_i}^2}\right)\left(w_i \frac{\mA_i^\trans\mI_i}{\norm{\mA_i}^2}\right)} + (1 - \frac{1}{\nthreads}) \E{ w_i \frac{\mI_i^\trans\mA_i}{\norm{\mA_i}^2}}\E{w_j \frac{\mA_j^\trans\mI_j}{\norm{\mA_j}^2}} \\
    & = \frac{1}{\nthreads} \E{w_i^2 \frac{\mI_i^\trans\mI_i}{\norm{\mA_i}^2}} + \left(1 - \frac{1}{\nthreads}\right) \mP  \mW\mD^{-2} \mA\mA^\trans \mP \mW \mD^{-2}\\
    & = \frac{1}{\nthreads} \mP \mW^2 \mD^{-2}+ \left(1 - \frac{1}{\nthreads}\right) \mP \mW\mD^{-2} \mA\mA^\trans \mP \mW \mD^{-2}
\end{align*}
by separating the cases where $i=j$ from those where $i \neq j$ and utilizing the independence of the indices sampled in $\tau_k$.

\section{Proof of \Cref{thm:expected_error}}\label{sec:expected_error_proof}
We prove \Cref{thm:expected_error} starting from from the error update in \Cref{eqn:error-update-matrices}. Expanding the squared error norm,
\begin{align*}
    \norm{e^{k+1}}^2 &= \norm{(\mI - \mA^\trans \mM_k \mA)e^{k} + \mA^\trans \mM_k r^\star}^2 \\
    &  =  \norm{(\mI - \mA^\trans \mM_k \mA)e^{k}}^2  + 2\ve{(\mI - \mA^\trans \mM_k \mA)e^{k}}{\mA^\trans \mM_k r^\star}  + \norm{\mA^\trans \mM_k r^\star}^2.
\end{align*}
Upon taking expections, the middle term simplifies since $\mA^\trans \E{\mM_k}r^\star=0$ by \Cref{ass:propI}. Thus,
\begin{equation}\label{eqn:before_taking_exp}
\begin{split}
    \E{\norm{e^{k+1}}^2}
    &= \E{\norm{(\mI - \mA^\trans \mM_k \mA)e^{k}}^2} - 2\E{\ve{\mA^\trans \mM_k \mA e^{k}}{\mA^\trans \mM_k r^\star}} + \E{\norm{\mA^\trans \mM_k r^\star}^2}.
\end{split}
\end{equation}

Making use of \Cref{lem:EMk} to take the expectation of the first term in \Cref{eqn:before_taking_exp},
\begin{align*}
    &\E{\norm{(\mI - \mA^\trans \mM_k \mA)e^{k}}^2 } \\
    &\quad= \E{\bigg \langle e^k, (\mI - \mA^\trans \mM_k \mA)^\trans (\mI - \mA^\trans \mM_k \mA)e^{k} \bigg \rangle}\\
    &\quad= \bigg \langle e^k, (\mI - 2\mA^\trans \E{\mM_k} \mA + \mA^\trans \E{\mM_k^\trans \mA \mA^\trans\mM_k} \mA)e^{k} \bigg \rangle\\
    &\quad= \bigg \langle e^k, \left(\mI - 2\alpha \frac{\mA^\trans\mA}{\norm{\mA}_F^2} + \frac{\alpha}{\nthreads}\frac{\mA^\trans\mW \mA}{\norm{\mA}_F^2} + \alpha^2\left(1 - \frac{1}{\nthreads}\right)  \left(\frac{\mA^\trans\mA}{\norm{\mA}_F^2}\right)^2 \right)e^{k}\bigg \rangle \\
    &\quad= \bigg \langle  e^k, \left(\left(\mI - \alpha \frac{\mA^\trans\mA}{\norm{\mA}_F^2}\right)^2 + \frac{\mA^\trans}{\norm{\mA}_F} \left(\frac{\alpha}{\nthreads} \mW - \frac{\alpha^2}{\nthreads} \frac{\mA\mA^\trans}{\norm{\mA}_F^2}\right) \frac{\mA}{\norm{\mA}_F} \right)e^{k}\bigg \rangle .
\end{align*}
Since $\mA^\trans r^\star = 0,$ for the second term,
\begin{align*}
    2\E{\ve{\mA^\trans \mM_k \mA e^{k}}{\mA^\trans \mM_k r^\star}} &= 2\ve{\mA e^{k}}{\E{\mM_k^\trans \mA \mA^\trans \mM_k }r^\star} \\
    &=  2\frac{\alpha}{\nthreads\norm{\mA}_F^2}\ve{\mA e^k}{ \mW r^\star}.
\end{align*}
Similarly, for the last term,
\begin{align*}
    \E{\norm{\mA^\trans \mM_k r^\star}^2}
    &= \frac{\alpha}{\nthreads}\frac{\norm{r^\star}_\mW^2}{\norm{\mA}_F^2} .
\end{align*}

Combining these in \Cref{eqn:before_taking_exp}, 
\begin{align*}
    \E{\norm{e^{k+1}}^2} &= \bve{e^k}{ \left(\mI - \alpha \frac{\mA^\trans\mA}{\norm{\mA}_F^2}\right)^2e^{k}} \\
    & \quad+\bve{e^k}{\frac{\mA^\trans}{\norm{\mA}_F^2} \left(\frac{\alpha}{\nthreads} \mW - \frac{\alpha^2}{\nthreads} \frac{\mA\mA^\trans}{\norm{\mA}_F^2}\right) \mA e^{k}}
    -2\frac{\alpha}{\nthreads} \frac{\ve{\mA e^k}{ \mW r^\star}}{\norm{\mA}_F^2} + \frac{\alpha}{\nthreads}\frac{\norm{r^\star}_\mW^2}{\norm{\mA}_F^2} \\
    &= \bve{e^k}{\left( \left(\mI - \alpha \frac{\mA^\trans\mA}{\norm{\mA}_F^2}\right)^2 -  \frac{\alpha^2}{\nthreads} \left(\frac{\mA^\trans \mA}{\norm{\mA}_F^2}\right)^2\right)e^{k}} 
     + \frac{\alpha}{\nthreads}\frac{\norm{r^k}_\mW^2}{\norm{\mA}_F^2} \\
    &\le \sigmax{\left(\mI - \alpha \frac{\mA^\trans\mA}{\norm{\mA}_F^2}\right)^2 -  \frac{\alpha^2}{\nthreads} \left(\frac{\mA^\trans\mA}{\norm{\mA}_F^2}\right)^2}\norm{e^{k}}^2
     + \frac{\alpha}{\nthreads}\frac{ \norm{r^k}_\mW^2}{\norm{\mA}_F^2}.
\end{align*}

\section{Proof of \Cref{thm:unif_opt_alpha}}
\label{sec:unif_opt_alpha_proof}

\begin{proof}
We seek to optimize the convergence rate constant from \Cref{cor:unif_weight_err},
\begin{align*}
    \sigma_{\max}\left(\left(\mI - \alpha \frac{\mA^\trans\mA}{\norm{\mA}_F^2}\right)^2 
    + \frac{\alpha^2 }{\nthreads}\left( \mI -  \frac{\mA^\trans\mA}{\norm{\mA}_F^2}\right)\frac{\mA^\trans\mA}{\norm{\mA}_F^2}\right)
\end{align*}
with respect to $\alpha$. To do this, we first simplify from  a matrix polynomial to a maximum over scalar polynomials in $\alpha$ with coefficients based on each singular value of $\mA$. We then show that the maximum occurs when either the minimum or maximum singular value of $\mA$ is used. Finally, we derive a condition for which singular value to use, and determine the optimal $\alpha$ that minimizes the maximum singular value.

Defining $\mQ^\trans \Sigma \mQ = \frac{\mA^\trans\mA}{\norm{\mA}_F^2}$ as the eigendecomposition, and the polynomial
\begin{equation*}
p(\sigma) \eqdef 1-2\alpha \sigma+\alpha^2\left(\frac{\sigma}{\nthreads} +\left(1-\frac{1}{\nthreads}\right)\sigma^2\right),
\end{equation*}
the convergence rate constant from \Cref{cor:unif_weight_err} can be written as $\sigmax{p\left(\frac{\mA^\trans\mA}{\norm{\mA}_F^2}\right)}$. Since $p\left(\frac{\mA^\trans\mA}{\norm{\mA}_F^2}\right)$ is a polynomial of a symmetric matrix, its singular vectors are the same as those of its argument, while its corresponding singular values are the polynomial $p$ applied to the singular values of the original matrix. That is,
\begin{equation*}
    p\left(\frac{\mA^\trans\mA}{\norm{\mA}_F^2}\right) = p\left(\mQ^\trans \Sigma \mQ\right) = \mQ^\trans p\left( \Sigma \right) \mQ.
\end{equation*}
Thus, the convergence rate constant can be written as
\begin{equation*}\label{eqn:singular-mat-polynomial}
    \sigmax{p\left(\frac{\mA^\trans\mA}{\norm{\mA}_F^2}\right)} = \sigmax{p(\Sigma)}.
\end{equation*}

Moreover, we can bound this extremal singular value by the maximum of the polynomial $p$ over an interval containing the spectrum of $\Sigma$
\begin{align*}
    \sigmax{p\left(\Sigma\right)} \leq \max \abs{p\left(\sigma\right)} \quad \text{subject to} \quad \sigma \in \left[s_{\min}, s_{\max}\right].
\end{align*}
Here, the singular values of $\Sigma$ are bounded from below by $s_{\min}\eqdef\frac{\sigminsqA}{\norm{\mA}_F^2}$ and above by $s_{\max}\eqdef\frac{\sigmaxsqA}{\norm{\mA}_F^2}$ since $\Sigma$ is the diagonal matrix of singular values of $\frac{\mA^\trans\mA}{\norm{\mA}_F^2}$.
Note that the polynomial can be factored as $p(\sigma) = (1-\sigma\alpha)^2+\frac{\sigma\alpha^2}{\nthreads}\left(1-\sigma\right)$, and is positive for $\sigma \in \left[0, 1\right]$, which contains $\left[s_{\min}, s_{\max}\right]$.
Also, since the coefficient of the $\sigma^2$ term of the polynomial $p$ is $\alpha^2 \left( 1 - \frac{1}{\nthreads} \right)$ which is greater than or equal to zero, the polynomial is convex in $\sigma$ on the interval $\left[s_{\min}, s_{\max}\right]$. Thus, the maximum of $p$ on the interval $\left[s_{\min}, s_{\max}\right]$ is attained at one of the two endpoints $s_{\min}, s_{\max}$ and we have the bound
\begin{align*}
    \sigma_{\max}\left(p\left(\Sigma\right)\right) = \max \left( p\left(s_{\min}\right), p\left(s_{\max}\right)\right).
\end{align*}

To optimize this bound with respect to $\alpha$, we first find conditions on $\alpha$ such that $p(s_{\min}) < p(s_{\max})$. If $s_{max}=s_{min}$, this obviously never holds; otherwise, $s_{max}>s_{min}$ and
\begin{align*}
    p(s_{\min}) &< p(s_{\max})\\
    1-2\alpha s_{\min}+\alpha^2\left[\frac{s_{\min}}{\nthreads} +\left(1-\frac{1}{\nthreads}\right)s_{\min}^2\right] &<
    1-2\alpha s_{\max}+\alpha^2\left[\frac{s_{\max}}{\nthreads} +\left(1-\frac{1}{\nthreads}\right)s_{\max}^2\right]
\end{align*}
Grouping like terms and cancelling, we get
\begin{align*}
    \alpha\left(2 - \frac{\alpha}{\nthreads} \right) \left(s_{\max} - s_{\min}\right) &< 
    \alpha^2\left(1 - \frac{1}{\nthreads}\right) \left(s_{\max}^2 - s_{\min}^2\right)
\end{align*}
Since $\frac{\alpha}{\nthreads}>0$, we can divide it from both sides.
\begin{align*}
    \left(2\nthreads  - \alpha \right) \left(s_{\max} - s_{\min}\right) &< 
    \alpha\left(\nthreads - 1\right) \left(s_{\max}^2 - s_{\min}^2\right)
\end{align*}
Since $s_{max} > s_{min}$, we can divide both sides by $s_{max}-s_{min}$.
\begin{align*}
    2 \nthreads - \alpha &< \alpha\left(\nthreads-1\right)\left(s_{\max}+s_{\min}\right) \\
    2q &< \alpha \left(1+\left(\nthreads-1\right)\left(s_{\max}+s_{\min}\right)\right)\\
    \alpha &> \frac{2 \nthreads}{1 + \left(\nthreads-1\right) \left(s_{\min} + s_{\max}\right)} \eqdef \widehat{\alpha}
\end{align*}

Thus,

\begin{align*}
    \sigma_{\max}\left(p\left(\Sigma\right)\right) = \begin{cases}
        p\left(s_{\max}\right), &\quad\alpha \geq \widehat{\alpha} \\
        p\left(s_{\min}\right), &\quad\alpha<\widehat{\alpha}
    \end{cases}
\end{align*}

For the first term,
\begin{align*}
    \frac{\partial}{\partial\alpha} p(s_{\max}) &= -2s_{\max} + 2 \left(\frac{s_{\max}}{\nthreads} + \left(1-\frac{1}{\nthreads}\right)s_{\max}^2\right)\alpha \\
    &\geq -2s_{\max} + 2\left(\frac{s_{\max}}{\nthreads} + \left(1-\frac{1}{\nthreads}\right)s_{\max}^2\right) \widehat{\alpha}
\end{align*}
since $\alpha \leq \widehat{\alpha}$ and the coefficient is positive. Factoring $\frac{2 s_{max}}{\nthreads}$ from the second term and substituting for $\widehat{\alpha}$, we get
\begin{align*}
    &= -2s_{\max} + \frac{2s_{\max}}{\nthreads}\left(1+ \left(\nthreads-1\right)s_{\max}\right)\widehat{\alpha} \\
    &= -2s_{\max} + \frac{2s_{\max}}{\nthreads}\left(1+ \left(\nthreads-1\right)s_{\max}\right)\frac{2 \nthreads}{1 + \left(1 - \nthreads\right) \left(s_{\min} + s_{\max}\right)} \\
    &= -2 s_{\max} + 2 s_{\max}\frac{2 \left(1+(\nthreads-1)s_{\max}\right)}{1+(\nthreads-1)(s_{\max}+s_{\min})} \\
    &= 2 s_{\max} \left[ -1 + \frac{2 \left(1+(\nthreads-1)s_{\max}\right)}{1+(\nthreads-1)(s_{\max}+s_{\min})}\right] \\
    &= 2 s_{\max} \left[\frac{ 1+(\nthreads-1)(s_{\max}-s_{\min})}{1+(\nthreads-1)(s_{\max}+s_{\min})}\right]\\
    &> 0
\end{align*}
since all terms in both numerator and denominator are positive. Thus, the function is monotonic increasing on $\alpha \in [\widehat{\alpha},\infty)$, and the minimum is at the lower endpoint, i.e. $\alpha^\star = \widehat{\alpha}$. \\
Similarly, for the second term,
\begin{align*}
    \frac{\partial}{\partial\alpha} p(s_{\min}) &= -2s_{\min} + 2 \left(\frac{s_{\min}}{\nthreads} + \left(1-\frac{1}{\nthreads}\right)s_{\min}^2\right)\alpha \\
    &< -2s_{\min} + 2\left(\frac{s_{\min}}{\nthreads} + \left(1-\frac{1}{\nthreads}\right)s_{\min}^2\right) \widehat{\alpha} \\
    &= 2 s_{\min} \left[\frac{ 1-(\nthreads-1)(s_{\max}-s_{\min})}{1+(\nthreads-1)(s_{\max}+s_{\min})}\right]
\end{align*}
If 
\begin{align} \label{opt_alpha_cond}
    1-(\nthreads-1)(s_{\max}-s_{\min})<0,
\end{align}
this function is monotonic decreasing on $\alpha \in (-\infty,\alpha^\star]$, and the minimum is at the upper endpoint i.e. $\alpha = \alpha^\star$. Otherwise, the minimum occurs at the critical point, so we set the derivative to 0 and solve for $\alpha^\star$
\begin{align*}
    \frac{\partial}{\partial\alpha} p(s_{\min}) &= -2s_{\min} + 2 \left(\frac{s_{\min}}{\nthreads} + \left(1-\frac{1}{\nthreads}\right)s_{\min}^2\right)\alpha^\star \\
    &= -2s_{\min} + \frac{2s_{\min}}{\nthreads}\left(1+ \left(\nthreads-1\right)s_{\min}\right)\alpha^\star \\
    &= 0 \\
    \frac{2s_{\min}}{\nthreads}\left(1 + \left(\nthreads-1\right)s_{\min}\right) \alpha^\star &= 2s_{\min} \\
    \alpha^\star &= \frac{\nthreads}{1+(\nthreads-1)s_{\min}}
\end{align*}
\end{proof}

\section{Corollary Proofs}
We provide proofs for the corollaries of \Cref{sec:conv}, which follow from \Cref{thm:expected_error}.

\subsection{Proof of \Cref{cor:unif_weight_err}}\label{subsec:unif_weight_cor_proof}

Suppose $p_i = \frac{\norm{\mA_i}^2}{\norm{\mA}_F^2}$ and $\mW = \alpha\mI$. 
From the proof of \Cref{thm:expected_error},
\begin{align*}
    \E{\norm{e^{k+1}}^2} 
    &= \bve{e^k}{\left( \left(\mI - \alpha \frac{\mA^\trans\mA}{\norm{\mA}_F^2}\right)^2 -  \frac{\alpha^2}{\nthreads} \left(\frac{\mA^\trans \mA}{\norm{\mA}_F^2}\right)^2\right)e^{k}} 
    + \frac{\alpha}{\nthreads}\frac{\norm{r^k}_\mW^2}{\norm{\mA}_F^2}.
\end{align*}
In this case, since $\mA^\trans r^\star = 0$, $\ve{\mA e^k}{r^\star} = 0$ and
\begin{align*}
    \norm{r^k}_\mW^2 &= \alpha\norm{\mA e^k}^2 + 2 \alpha\ve{\mA e^k}{r^\star} + \alpha\norm{r^\star}^2 \\
    &= \alpha\ve{e^k}{\mA^\trans \mA e^k}+ \alpha\norm{r^\star}^2.
\end{align*}
Combining the inner products, 
\begin{align*}
    \E{\norm{e^{k+1}}^2}
    &= \bigg \langle  e^k,\left( \left(\mI - \alpha \frac{\mA^\trans\mA}{\norm{\mA}_F^2}\right)^2 
    + \frac{\alpha^2}{\nthreads} \left(\mI - \frac{\mA^\trans \mA}{\norm{\mA}_F^2}\right)\frac{\mA^\trans \mA}{\norm{\mA}_F^2}\right)e^{k}\bigg  \rangle
    + \frac{\alpha^2\norm{r^\star}^2}{\nthreads\norm{\mA}_F^2}\\
    &\le \sigma_{\max}\left(\left(\mI - \alpha \frac{\mA^\trans\mA}{\norm{\mA}_F^2}\right)^2 
    + \frac{\alpha^2 }{\nthreads}\left( \mI -  \frac{\mA^\trans\mA}{\norm{\mA}_F^2}\right)\frac{\mA^\trans\mA}{\norm{\mA}_F^2}\right)\norm{e^{k}}^2  + \frac{\alpha^2\norm{r^\star}^2}{\nthreads\norm{\mA}_F^2} .
\end{align*}

\subsection{Proof of \Cref{cor:RK_conv}}\label{subsec:RK_conv_proof}
Suppose $\nthreads = 1$, $\mW=\mI$ and $p_i = \frac{\norm{\mA_i}^2}{\norm{\mA}_F^2}$. 
\begin{align*}
    \E{\norm{e^{k+1}}^2} 
    &\le \sigmax{\mI -  \frac{\mA^\trans\mA}{\norm{\mA}_F^2}}\norm{e^{k}}^2  + \frac{ \norm{r^\star}^2}{\norm{\mA}_F^2}
     \\ &= \left(1 -  \frac{\sigminsqA}{\norm{\mA}_F^2}\right)\norm{e^{k}}^2  + \frac{ \norm{r^\star}^2}{\norm{\mA}_F^2}
    .
\end{align*}
From the proof of \Cref{thm:expected_error},
\begin{align*}
    \E{\norm{e^{k+1}}^2}
    &= \bve{e^k}{\left( \left(\mI -  \frac{\mA^\trans\mA}{\norm{\mA}_F^2}\right)^2 -   \left(\frac{\mA^\trans \mA}{\norm{\mA}_F^2}\right)^2\right)e^{k}}
     + \frac{\norm{r^k}^2}{\norm{\mA}_F^2}.
\end{align*}
Decomposing $r^k$, 
\begin{align*}
    \norm{r^k}^2 &= \norm{\mA e^k}^2 + \norm{r^\star}^2 \\ &=\ve{e^k}{\mA^\trans\mA e^k} + \norm{r^\star}^2 .
\end{align*}
Combining the inner products,
\begin{align*}
    \E{\norm{e^{k+1}}^2}
    &= \bigg \langle e^k, \left( \left(\mI - \frac{\mA^\trans\mA}{\norm{\mA}_F^2}\right)^2 
    -\left(\frac{\mA^\trans \mA}{\norm{\mA}_F^2} \right)^2+ \frac{\mA^\trans\mA}{\norm{\mA}_F^2}\right)e^{k}\bigg \rangle + \frac{\norm{r^\star}^2}{\norm{\mA}_F^2} \\
    &= \bve{e^k}{ \left(\mI - \frac{\mA^\trans\mA}{\norm{\mA}_F^2}\right)e^{k}} + \frac{\norm{r^\star}^2}{\norm{\mA}_F^2}\\
     &\le \sigmax{\mI -  \frac{\mA^\trans\mA}{\norm{\mA}_F^2}}\norm{e^{k}}^2  + \frac{ \norm{r^\star}^2}{\norm{\mA}_F^2}
     \\ 
     &= \left(1 -  \frac{\sigminsqA}{\norm{\mA}_F^2}\right)\norm{e^{k}}^2  + \frac{ \norm{r^\star}^2}{\norm{\mA}_F^2}
    .
\end{align*}
\end{document}